\documentclass[12pt,reqno]{amsart}

\newcommand{\df}{\dfrac}
\newcommand{\tf}{\tfrac}

 \renewcommand{\a}{\alpha}
\renewcommand{\b}{\beta}

\renewcommand{\(}{\left\(}
\renewcommand{\)}{\right\)}
\renewcommand{\[}{\left\[}
\renewcommand{\]}{\right\]}

\numberwithin{equation}{section}
 \theoremstyle{plain}
\newtheorem{theorem}{Theorem}[section]

\newtheorem{corollary}[theorem]{Corollary}

\begin{document}

\title{ANALOGUES OF A TRANSFORMATION FORMULA OF RAMANUJAN}
\author{Atul Dixit}
\address{Department of Mathematics, University of Illinois, 1409 West Green
Street, Urbana, IL 61801, USA} \email{aadixit2@illinois.edu}
\maketitle
\centerline{\emph{Dedicated to Professor Bruce C.~Berndt on the occasion of his 70th birthday}}
\begin{abstract}
We derive two new analogues of a transformation formula of Ramanujan involving the Gamma and Riemann zeta functions present in the Lost Notebook. Both involve infinite series consisting of Hurwitz zeta functions and yield modular relations. As a special case of the first formula, we obtain an identity involving polygamma functions given by A.P.~Guinand and as a limiting case of the second formula, we derive the transformation formula of Ramanujan.
\end{abstract}
\section{Introduction}
In the volume \cite{lnb} containing Ramanujan's Lost Notebook are present some manuscripts of Ramanujan in the handwriting of G.N.~Watson. The first of these manuscripts contains the following beautiful claim (see \cite[p.~220]{lnb}).
\begin{theorem}\label{entry1} Define
\begin{equation}\label{w1.27}
\phi(x):=\psi(x)+\df{1}{2x}-\log x,
\end{equation}
where
\begin{equation}\label{w1.15b}
\psi(x):=\df{\Gamma^\prime(x)}{\Gamma(x)}=-\gamma-\sum_{m=0}^{\infty}\left(\df{1}{m+x}-\df{1}{m+1}\right),
\end{equation}
the logarithmic derivative of Gamma function.
Let the Riemann's $\xi$-function be defined by
\begin{equation*}
\xi(s):=(s-1)\pi^{-\tf{1}{2}s}\Gamma(1+\tf{1}{2}s)\zeta(s),
\end{equation*}
and let 
\begin{equation}\label{xif}
\Xi(t):=\xi(\tf{1}{2}+it)
\end{equation}
be the Riemann $\Xi$-function.
If $\a$ and $\b$ are positive numbers such that $\a\b=1$, then
\begin{multline}\label{w1.26}
\sqrt{\a}\left\{\df{\gamma-\log(2\pi\a)}{2\a}+\sum_{n=1}^{\infty}\phi(n\a)\right\}
=\sqrt{\b}\left\{\df{\gamma-\log(2\pi\b)}{2\b}+\sum_{n=1}^{\infty}\phi(n\b)\right\}\\
=-\df{1}{\pi^{3/2}}\int_0^{\infty}\left|\Xi\left(\df{1}{2}t\right)\Gamma\left(\df{-1+it}{4}\right)\right|^2
\df{\cos\left(\tf{1}{2}t\log\a\right)}{1+t^2}\, dt,
\end{multline}
where $\gamma$ denotes Euler's constant.
\end{theorem}
A.P.~Guinand \cite{apg2, apg3} rediscovered the first equality in (\ref{w1.26}) in a slightly different form. Recently, B.C.~Berndt and A.~Dixit \cite{bcbad} proved both parts of (\ref{w1.26}). A key element in their proof was the identity \cite[p.~260, eqn.~(22)]{riemann} or \cite[p.~77, eqn.~(22)]{cp}, true for $n$ real,
\begin{align}\label{Ramint}
&\int_{0}^{\infty}\Gamma\left(\frac{-1+it}{4}\right)\Gamma\left(\frac{-1-it}{4}\right)
\left(\Xi\left(\frac{1}{2}t\right)\right)^2\frac{\cos nt}{1+t^2}\, dt\nonumber\\
&\hspace{1cm}=\pi^{3/2}\int_{0}^{\infty}\left(\frac{1}{e^{xe^{n}}-1}-\frac{1}{xe^{n}}\right)
\left(\frac{1}{e^{xe^{-n}}-1}-\frac{1}{xe^{-n}}\right)\, dx.
\end{align}
Ramanujan's paper \cite{riemann} consists of other equations similar to (\ref{Ramint}). Motivated by the use of (\ref{Ramint}) in deriving (\ref{w1.26}), we work with two other equations, namely equations (19) and (20) in \cite{riemann}, to derive two new analogues of (\ref{w1.26}). However, it must be pointed out that equation (19) in \cite{riemann}, as it stands, is incorrect. The second term on its right-hand side, namely, $\displaystyle-\frac{1}{4}(4\pi)^{\frac{(s-3)}{2}}\Gamma(s)\zeta(s)\cosh n(1-s)$ should not be present. Equations (19) (corrected) and (20) are respectively as follows:

\begin{theorem}
For Re $s>1$ and $n$ real, we have
\begin{align}\label{eq19cor}
&\int_{0}^{\infty}\Gamma\left(\frac{s-1+it}{4}\right)\Gamma\left(\frac{s-1-it}{4}\right)
\Xi\left(\frac{t+is}{2}\right)\Xi\left(\frac{t-is}{2}\right)\frac{\cos nt}{(s+1)^2+t^2}\, dt\nonumber\\
&=\frac{1}{8}(4\pi)^{-\frac{(s-3)}{2}}\int_{0}^{\infty}\frac{x^{s}}{(e^{xe^{n}}-1)(e^{xe^{-n}}-1)}\, dx,
\end{align}
where $\Xi(t)$ is as defined in (\ref{xif}).
\end{theorem}

\begin{theorem}
For $-1<$ Re $s<1$ and $n$ real, we have
\begin{align}\label{eq20}
&\int_{0}^{\infty}\Gamma\left(\frac{s-1+it}{4}\right)\Gamma\left(\frac{s-1-it}{4}\right)
\Xi\left(\frac{t+is}{2}\right)\Xi\left(\frac{t-is}{2}\right)\frac{\cos nt}{(s+1)^2+t^2}\, dt\nonumber\\
&=\frac{1}{8}(4\pi)^{-\frac{(s-3)}{2}}\int_{0}^{\infty}x^{s}\left(\frac{1}{e^{xe^{n}}-1}-\frac{1}{xe^{n}}\right)\left(\frac{1}{e^{xe^{-n}}-1}-\frac{1}{xe^{-n}}\right)\, dx,
\end{align}
where again $\Xi(t)$ is as defined in (\ref{xif}).
\end{theorem}
The identity (\ref{Ramint}) is the special case $s=0$ of (\ref{eq20}). Now we state the two key theorems in this paper which give two new analogues of (\ref{w1.26}).

\begin{theorem}\label{mainn}
Let $\zeta(z,a)$ denote the Hurwitz zeta function defined for Re $z>1$ by
\begin{equation}
\zeta(z,a)=\sum_{n=0}^{\infty}\frac{1}{(n+a)^z}.
\end{equation}
If $\a$ and $\b$ are positive numbers such that $\a\b=1$, then for Re $z>2$ and $1<c<$ Re $z-1$,
\begin{align}\label{maineqnr}
&\a^{-\frac{z}{2}}\sum_{k=1}^{\infty}\zeta\left(z,1+\frac{k}{\a}\right)=\b^{-\frac{z}{2}}\sum_{k=1}^{\infty}\zeta\left(z,1+\frac{k}{\b}\right)\nonumber\\
&=\frac{\alpha^{\frac{z}{2}}}{2\pi i\Gamma(z)}\int_{c-i\infty}^{c+i\infty}\Gamma(s)\zeta(s)\Gamma(z-s)\zeta(z-s)\a^{-s}\, ds\nonumber\\
&=\frac{8(4\pi)^{\frac{z-4}{2}}}{\Gamma(z)}\int_{0}^{\infty}\Gamma\left(\frac{z-2+it}{4}\right)\Gamma\left(\frac{z-2-it}{4}\right)
\Xi\left(\frac{t+i(z-1)}{2}\right)\Xi\left(\frac{t-i(z-1)}{2}\right)\frac{\cos\left( \tf{1}{2}t\log\a\right)}{z^2+t^2}\, dt,
\end{align}
where $\Xi$ is defined as in (\ref{xif}).
\end{theorem}

\begin{theorem}
Let $0<$ Re $z<2$. Define $\varphi(z,x)$ as
\begin{equation}
\varphi(z,x)=\zeta(z,x)-\frac{1}{2}x^{-z}+\frac{x^{1-z}}{1-z},
\end{equation}
where $\zeta(z,x)$ denotes the Hurwitz zeta function. Then if $\alpha$ and $\beta$ are any positive numbers such that $\alpha\beta=1$,
\begin{align}\label{mainneq2}
&\a^{\frac{z}{2}}\left(\sum_{n=1}^{\infty}\varphi(z,n\a)-\frac{\zeta(z)}{2\alpha^{z}}-\frac{\zeta(z-1)}{\alpha (z-1)}\right)=\b^{\frac{z}{2}}\left(\sum_{n=1}^{\infty}\varphi(z,n\b)-\frac{\zeta(z)}{2\beta^{z}}-\frac{\zeta(z-1)}{\beta (z-1)}\right)\nonumber\\
&=\frac{8(4\pi)^{\frac{z-4}{2}}}{\Gamma(z)}\int_{0}^{\infty}\Gamma\left(\frac{z-2+it}{4}\right)\Gamma\left(\frac{z-2-it}{4}\right)
\Xi\left(\frac{t+i(z-1)}{2}\right)\Xi\left(\frac{t-i(z-1)}{2}\right)\frac{\cos\left( \tf{1}{2}t\log\a\right)}{z^2+t^2}\, dt,
\end{align}
where $\Xi(t)$ is defined in (\ref{xif}).
\end{theorem}

This paper is organized as follows. In Section 2, we review some basic properties of Mellin transforms which are used subsequently. Then in Section 3, we derive an analogue of (\ref{w1.26}), namely (\ref{maineqnr}), using two different methods, both of which make use of (\ref{eq19cor}). In Section 4, we derive a second analogue of (\ref{w1.26}), namely (\ref{mainneq2}), which makes use of (\ref{eq20}), and gives (\ref{w1.26}) as a limiting case. Finally in Section 5, we prove (\ref{eq19cor}). 

\section{Basic Properties of Mellin Transforms}
Let $F(z)$ denote the Mellin transform of $f(x)$, i.e.,
\begin{equation}\label{mt}
F(z)=\int_{0}^{\infty}x^{z-1}f(x)\, dx.
\end{equation}
Then the inverse Mellin transform is given by
\begin{equation}\label{imt}
f(x)=\frac{1}{2\pi i}\int_{c-i\infty}^{c+i\infty}F(z)x^{-z}\, dz,
\end{equation}
where $c$ lies in the fundamental strip (or the strip of analyticity) for which $F(z)$ is defined.
We also note the Mellin convolution theorem \cite[p.~83]{kp} which states that if $F(z)$ and $G(z)$ are Mellin transforms of $f(x)$ and $g(x)$ respectively, then
\begin{equation}\label{mct}
\int_{0}^{\infty}x^{z-1}f(x)g(x)\, dx=\frac{1}{2\pi i}\int_{c-i\infty}^{c+i\infty}F(s)G(z-s)\, ds,
\end{equation}
$c$ again being in the associated fundamental strip.

Now let $F(z)$ be related to $f(x)$ by (\ref{mt}) and (\ref{imt}), where $f(x)$ is locally integrable on $(0,\infty)$, is $O(x^{-a})$ as $x\to 0^{+}$ and $O(x^{-b})$, where $b>1$, as $x\to\infty$ and $a<c<b$. Then for max$\{1,a\}<c<b$, we have
\begin{equation}\label{melser}
\sum_{n=1}^{\infty}f(nx)=\frac{1}{2\pi i}\int_{c-i\infty}^{c+i\infty}F(s)\zeta(s)x^{-s}\, ds,
\end{equation}
where $\zeta(s)$ denotes the Riemann zeta function (see \cite[p.~117]{kp}).

\section{The first modular relation involving Hurwitz zeta functions}
\begin{theorem}\label{main}
Let $\zeta(z,a)$ denote the Hurwitz zeta function defined for Re $z>1$ by
\begin{equation}
\zeta(z,a)=\sum_{n=0}^{\infty}\frac{1}{(n+a)^z}.
\end{equation}
If $\a$ and $\b$ are positive numbers such that $\a\b=1$, then for Re $z>2$ and $1<c<$ Re $z-1$,
\begin{align}\label{maineq}
&\a^{-\frac{z}{2}}\sum_{k=1}^{\infty}\zeta\left(z,1+\frac{k}{\a}\right)=\b^{-\frac{z}{2}}\sum_{k=1}^{\infty}\zeta\left(z,1+\frac{k}{\b}\right)\nonumber\\
&=\frac{\alpha^{\frac{z}{2}}}{2\pi i\Gamma(z)}\int_{c-i\infty}^{c+i\infty}\Gamma(s)\zeta(s)\Gamma(z-s)\zeta(z-s)\a^{-s}\, ds\nonumber\\
&=\frac{8(4\pi)^{\frac{z-4}{2}}}{\Gamma(z)}\int_{0}^{\infty}\Gamma\left(\frac{z-2+it}{4}\right)\Gamma\left(\frac{z-2-it}{4}\right)
\Xi\left(\frac{t+i(z-1)}{2}\right)\Xi\left(\frac{t-i(z-1)}{2}\right)\frac{\cos\left( \tf{1}{2}t\log\a\right)}{z^2+t^2}\, dt,
\end{align}
where $\Xi$ is defined as in (\ref{xif}).
\end{theorem}
\textbf{First proof:}
\vspace{2mm}
Replace $s$ by $z-1$ in (\ref{eq19cor}) and then multiply the resulting two sides by $\displaystyle 8(4\pi)^{\frac{z-4}{2}}e^{-nz}$, so that for Re $z>2$ and $n$ real, we have
\begin{align}\label{eq19cor2sec}
&8(4\pi)^{\frac{z-4}{2}}e^{-nz}\int_{0}^{\infty}\Gamma\left(\frac{z-2+it}{4}\right)\Gamma\left(\frac{z-2-it}{4}\right)
\Xi\left(\frac{t+i(z-1)}{2}\right)\Xi\left(\frac{t-i(z-1)}{2}\right)\frac{\cos nt}{z^2+t^2}\, dt\nonumber\\
&=e^{-nz}\int_{0}^{\infty}\frac{x^{z-1}}{(e^{xe^{n}}-1)(e^{xe^{-n}}-1)}\, dx.
\end{align}
The integral on the right-hand side of (\ref{eq19cor2sec}) can be viewed as the Mellin transform of the product of $\displaystyle\frac{1}{e^{xe^n}-1}$ and $\displaystyle\frac{1}{e^{xe^{-n}}-1}$.

But from \cite[p.~18, eqn. 2.4.1]{titch}, we know that for Re $z>1$,
\begin{equation}\label{gamzeta}
\Gamma(z)\zeta(z)=\int_{0}^{\infty}\frac{t^{z-1}}{e^{t}-1}\, dt.
\end{equation}
Employing a change of variable $t=xe^{n}$ in (\ref{gamzeta}), we deduce that
\begin{equation}\label{gamzeta2}
e^{-nz}\Gamma(z)\zeta(z)=\int_{0}^{\infty}\frac{x^{z-1}}{e^{xe^{n}}-1}\, dx.
\end{equation}
Similarly letting $t=xe^{-n}$ in (\ref{gamzeta}), we find that
\begin{equation}\label{gamzeta3}
e^{nz}\Gamma(z)\zeta(z)=\int_{0}^{\infty}\frac{x^{z-1}}{e^{xe^{-n}}-1}\, dx.
\end{equation}
Thus from (\ref{mct}), (\ref{gamzeta2}) and (\ref{gamzeta3}), we see that for $1<c<$ Re $z-1$,
\begin{equation}\label{melrelat}
\int_{0}^{\infty}\frac{x^{z-1}}{(e^{xe^{n}}-1)(e^{xe^{-n}}-1)}\, dx=\frac{1}{2\pi i}\int_{c-i\infty}^{c+i\infty}e^{-ns}\Gamma(s)\zeta(s)e^{n(z-s)}\Gamma(z-s)\zeta(z-s)\, ds,
\end{equation}
which can be written as
\begin{equation}\label{melrelat2}
e^{-nz}\int_{0}^{\infty}\frac{x^{z-1}}{(e^{xe^{n}}-1)(e^{xe^{-n}}-1)}\, dx=\frac{1}{2\pi i}\int_{c-i\infty}^{c+i\infty}e^{-2ns}\Gamma(s)\zeta(s)\Gamma(z-s)\zeta(z-s)\, ds,
\end{equation}
Letting $n=\tf{1}{2}\log\a$ in (\ref{eq19cor2sec}) and (\ref{melrelat2}) and combining them together, we arrive at
\begin{align}\label{eq19cor2com}
&8(4\pi)^{\frac{z-4}{2}}\a^{-\frac{z}{2}}\int_{0}^{\infty}\Gamma\left(\frac{z-2+it}{4}\right)\Gamma\left(\frac{z-2-it}{4}\right)
\Xi\left(\frac{t+i(z-1)}{2}\right)\Xi\left(\frac{t-i(z-1)}{2}\right)\frac{\cos\left(\tf{1}{2}t\log\a\right)}{z^2+t^2}\, dt\nonumber\\
&=\frac{1}{2\pi i}\int_{c-i\infty}^{c+i\infty}\Gamma(s)\zeta(s)\Gamma(z-s)\zeta(z-s)\a^{-s}\, ds.
\end{align}
Upon simplification, this gives the last equality in (\ref{maineq}).
Now since $1<c<$ Re $z-1$, on the vertical line Re $s=c$, we have Re $(z-s)>1$. So we can use the representation 
\begin{equation}\label{zeta}
\zeta(z-s)=\sum_{k=1}^{\infty}\frac{1}{k^{z-s}}.
\end{equation}

Using (\ref{zeta}) on the right-hand side of (\ref{eq19cor2com}), we find, by absolute convergence, that
\begin{align}\label{simpcon}
&\frac{1}{2\pi i}\int_{c-i\infty}^{c+i\infty}\Gamma(s)\zeta(s)\Gamma(z-s)\zeta(z-s)\a^{-s}\, ds=\sum_{k=1}^{\infty}k^{-z}\frac{1}{2\pi i}\int_{c-i\infty}^{c+i\infty}\Gamma(s)\zeta(s)\Gamma(z-s)\left(\frac{\a}{k}\right)^{-s}\, ds\nonumber\\
&=\Gamma(z)\sum_{k=1}^{\infty}k^{-z}\frac{1}{2\pi i}\int_{c-i\infty}^{c+i\infty}\left(\frac{\Gamma(s)\Gamma(z-s)}{\Gamma(z)}\right)\zeta(s)\left(\frac{\a}{k}\right)^{-s}\, ds.
\end{align}
Now we know that for $0<$ Re $s<$ Re $z$, Euler's beta integral $B(s,z-s)$ is given by
\begin{equation}\label{betamel}
B(s,z-s)=\int_{0}^{\infty}\frac{x^{s-1}}{(1+x)^{z}}\, dx=\frac{\Gamma(s)\Gamma(z-s)}{\Gamma(z)}.
\end{equation}
In other words, $B(s,z-s)$ is the Mellin transform of $\displaystyle\frac{1}{(1+x)^{z}}$.
Now for Re $z>2$, $f(x):=\displaystyle\frac{1}{(1+x)^{z}}$ is locally integrable on $(0,\infty)$. Also, as $x\to 0^{+}$, $f(x)=O(1)$ and as $x\to \infty$, $f(x)\sim \displaystyle\frac{1}{x^{z}}=O\left(\frac{1}{x^{\text{Re}(z)}}\right)$. In particular, for $1<c< $ Re $z-1$, using (\ref{melser}) and (\ref{betamel}), we find that
\begin{equation}\label{melser2}
\frac{1}{2\pi i}\int_{c-i\infty}^{c+i\infty}B(s,z-s)\zeta(s)x^{-s}\, ds=\sum_{m=1}^{\infty}(1+xm)^{-z}.
\end{equation}
The integral in the second expression in (\ref{simpcon}) can also be directly evaluated using formula 5.78 in \cite[p.~202]{ober}.\\

From (\ref{simpcon}), (\ref{betamel}) and (\ref{melser2}), we arrive at
\begin{align}\label{simpcon2}
&\frac{1}{2\pi i}\int_{c-i\infty}^{c+i\infty}\Gamma(s)\zeta(s)\Gamma(z-s)\zeta(z-s)\a^{-s}\, ds
=\Gamma(z)\sum_{k=1}^{\infty}k^{-z}\sum_{m=1}^{\infty}\left(1+\frac{\a m}{k}\right)^{-z}\nonumber\\
&=\Gamma(z)\sum_{k=1}^{\infty}\sum_{m=1}^{\infty}(k+\a m)^{-z}=\a^{-z}\Gamma(z)\sum_{k=1}^{\infty}\zeta\left(z,1+\frac{k}{\a}\right).
\end{align}
Invoking (\ref{simpcon2}) in (\ref{eq19cor2com}), simplifying and rearranging, we see that
\begin{align}\label{bf2}
&\a^{-\frac{z}{2}}\sum_{k=1}^{\infty}\zeta\left(z,1+\frac{k}{\a}\right)\nonumber\\
&=\frac{8(4\pi)^{\frac{(z-4)}{2}}}{\Gamma(z)}\int_{0}^{\infty}\Gamma\left(\frac{z-2+it}{4}\right)\Gamma\left(\frac{z-2-it}{4}\right)
\Xi\left(\frac{t+i(z-1)}{2}\right)\Xi\left(\frac{t-i(z-1)}{2}\right)\frac{\cos\left( \tf{1}{2}t\log\a\right)}{z^2+t^2}\, dt.
\end{align}
Now replacing $\a$ by $\b$ in (\ref{bf2}), we see that
\begin{align}\label{bf3}
&\b^{-\frac{z}{2}}\sum_{k=1}^{\infty}\zeta\left(z,1+\frac{k}{\b}\right)\nonumber\\
&=\frac{8(4\pi)^{\frac{(z-4)}{2}}}{\Gamma(z)}\int_{0}^{\infty}\Gamma\left(\frac{z-2+it}{4}\right)\Gamma\left(\frac{z-2-it}{4}\right)
\Xi\left(\frac{t+i(z-1)}{2}\right)\Xi\left(\frac{t-i(z-1)}{2}\right)\frac{\cos\left( \tf{1}{2}t\log\b\right)}{z^2+t^2}\, dt\nonumber\\
&=\frac{8(4\pi)^{\frac{(z-4)}{2}}}{\Gamma(z)}\int_{0}^{\infty}\Gamma\left(\frac{z-2+it}{4}\right)\Gamma\left(\frac{z-2-it}{4}\right)
\Xi\left(\frac{t+i(z-1)}{2}\right)\Xi\left(\frac{t-i(z-1)}{2}\right)\frac{\cos\left( \tf{1}{2}t\log\tf{1}{\a}\right)}{z^2+t^2}\, dt\nonumber\\
&=\frac{8(4\pi)^{\frac{(z-4)}{2}}}{\Gamma(z)}\int_{0}^{\infty}\Gamma\left(\frac{z-2+it}{4}\right)\Gamma\left(\frac{z-2-it}{4}\right)
\Xi\left(\frac{t+i(z-1)}{2}\right)\Xi\left(\frac{t-i(z-1)}{2}\right)\frac{\cos\left( \tf{1}{2}t\log\a\right)}{z^2+t^2}\, dt.\nonumber\\
\end{align}
Thus (\ref{bf2}) and (\ref{bf3}) establish (\ref{maineq}) and this completes the proof of Theorem \ref{main}.

\hfill $\square$

\textbf{Second proof:}
Letting $n=\tf{1}{2}\log\a$ in (\ref{eq19cor2sec}) and multiplying both sides by $\displaystyle\frac{1}{\Gamma(z)}$, we see that
\begin{align}\label{fin2}
&\frac{8(4\pi)^{\frac{(z-4)}{2}}}{\Gamma(z)}\int_{0}^{\infty}\Gamma\left(\frac{z-2+it}{4}\right)\Gamma\left(\frac{z-2-it}{4}\right)
\Xi\left(\frac{t+i(z-1)}{2}\right)\Xi\left(\frac{t-i(z-1)}{2}\right)\frac{\cos\left( \tf{1}{2}t\log\a\right)}{z^2+t^2}\, dt\nonumber\\
&=\frac{1}{\Gamma(z)}\int_{0}^{\infty}\frac{x^{z-1}}{(e^{x\sqrt{\a}}-1)(e^{x/\sqrt{\a}}-1)}\, dx\nonumber\\
&=\frac{\a^{-\frac{z}{2}}}{\Gamma(z)}
\int_{0}^{\infty}\frac{t^{z-1}}{(e^{t}-1)(e^{t/\a}-1)}\, dt,
\end{align}
where in the penultimate line, we have made a change of variable $t=x\sqrt{\a}$.

Now,
\begin{align}\label{left}
\frac{\a^{-\frac{z}{2}}}{\Gamma(z)}
\int_{0}^{\infty}\frac{t^{z-1}}{(e^{t}-1)(e^{t/\a}-1)}\, dt&=\frac{\a^{-\frac{z}{2}}}{\Gamma(z)}
\int_{0}^{\infty}\frac{t^{z-1}e^{-t}}{1-e^{-t}}\sum_{k=1}^{\infty}e^{-kt/\a}\, dt\nonumber\\
&=\frac{\a^{-\frac{z}{2}}}{\Gamma(z)}\sum_{k=1}^{\infty}
\int_{0}^{\infty}\frac{t^{z-1}e^{-(1+k/\a)t}}{1-e^{-t}}\, dt,
\end{align}
where the order of summation and integration can be interchanged because of absolute convergence.
But from \cite[p.~37, eqn. 2.17.1]{titch}, we know that for Re $z>1$,
\begin{equation}\label{hurwitzint}
\zeta(z,a)=\frac{1}{\Gamma(z)}\int_{0}^{\infty}\frac{x^{z-1}e^{-ax}}{1-e^{-x}}\, dx.
\end{equation}
Using (\ref{hurwitzint}) in (\ref{left}), we deduce that
\begin{equation}\label{used}
\frac{\a^{-\frac{z}{2}}}{\Gamma(z)}\int_{0}^{\infty}\frac{t^{z-1}}{(e^{t}-1)(e^{t/\a}-1)}\, dt=\a^{-\frac{z}{2}}\sum_{k=1}^{\infty}\zeta\left(z,1+\frac{k}{\a}\right).
\end{equation}
Thus from (\ref{fin2}) and (\ref{used}), we derive (\ref{bf2}). Following the same argument as shown in (\ref{bf3}), we obtain the other equality in (\ref{maineq}) as well. This finishes the second proof.

\hfill $\square$
\begin{corollary}\label{thm3.1spl}
For Re $z>2$, we have
\begin{align}\label{maineqspl}
&\sum_{k=1}^{\infty}\zeta\left(z,1+k\right)=\zeta(z-1)-\zeta(z)\nonumber\\
&=\frac{8(4\pi)^{\frac{z-4}{2}}}{\Gamma(z)}\int_{0}^{\infty}\Gamma\left(\frac{z-2+it}{4}\right)\Gamma\left(\frac{z-2-it}{4}\right)
\Xi\left(\frac{t+i(z-1)}{2}\right)\Xi\left(\frac{t-i(z-1)}{2}\right)\frac{dt}{z^2+t^2}.
\end{align}
\end{corollary}
\begin{proof}
Set $\a=1$ in (\ref{maineq}) and note that from \cite[p.~35]{titch}, for $1<c<$ Re $z-1$, we have 
\begin{equation}\label{titchmarsh}
\frac{1}{2\pi i}\int_{c-i\infty}^{c+i\infty}\Gamma(s)\zeta(s)\Gamma(z-s)\zeta(z-s)\, ds=\Gamma(z)\left(\zeta(z-1)-\zeta(z)\right).
\end{equation}
\end{proof}
\subsection{Guinand's formula as a special case of (\ref{maineq})}
Let $\psi^{(k)}(x)$ denote the $k^{th}$ derivative of the digamma function $\psi(x)$ defined in (\ref{w1.15b}), also known as the polygamma function of order $k$. In \cite{apg3}, Guinand gave the following formula
\begin{equation}\label{gui}
\sum_{n=1}^{\infty}\psi^{(k)}(1+nx)=x^{-k-1}\sum_{n=1}^{\infty}\psi^{(k)}\left(1+\frac{n}{x}\right),
\end{equation}
where $k\geq 2$. We derive this formula as a special case of (\ref{maineq}). Let $z$ be a natural number greater than $2$. From (\ref{w1.15b}), by successive differentiation, we find that
\begin{equation}\label{psidiff}
\psi^{(z-1)}(x)=(-1)^{z}(z-1)!\sum_{m=1}^{\infty}\frac{1}{(m-1+x)^{z}}.
\end{equation}
Thus,
\begin{align}\label{guileft}
\alpha^{\frac{z}{2}}\sum_{k=1}^{\infty}\psi^{(z-1)}(1+k\a)&=(-1)^{z}(z-1)!\alpha^{\frac{z}{2}}\sum_{k=1}^{\infty}\sum_{m=1}^{\infty}\frac{1}{(m+k\a)^{z}}\nonumber\\
&=(-1)^{z}(z-1)!\alpha^{\frac{z}{2}}\sum_{m=1}^{\infty}\sum_{k=1}^{\infty}\frac{1}{(m+k\a)^{z}}\nonumber\\
&=(-1)^{z}(z-1)!\alpha^{\frac{-z}{2}}\sum_{m=1}^{\infty}\sum_{k=1}^{\infty}\frac{1}{(k+m/\a)^{z}}\nonumber\\
&=(-1)^{z}(z-1)!\alpha^{\frac{-z}{2}}\sum_{m=1}^{\infty}\zeta\left(z,1+\frac{m}{\a}\right),
\end{align}
where the change in the order of summation in the second equality is justified by absolute convergence.

Then from (\ref{guileft}) and (\ref{maineq}), we obtain the following alternate version of (\ref{maineq}) when $z$ is a natural number greater than $2$,
\begin{align}\label{guigen}
&\alpha^{\frac{z}{2}}\sum_{k=1}^{\infty}\psi^{(z-1)}(1+k\a)=\b^{\frac{z}{2}}\sum_{k=1}^{\infty}\psi^{(z-1)}(1+k\b)\nonumber\\
&=8(-1)^z(4\pi)^{\frac{z-4}{2}}\nonumber\\
&\quad\times\int_{0}^{\infty}\Gamma\left(\frac{z-2+it}{4}\right)\Gamma\left(\frac{z-2-it}{4}\right)
\Xi\left(\frac{t+i(z-1)}{2}\right)\Xi\left(\frac{t-i(z-1)}{2}\right)\frac{\cos\left( \tf{1}{2}t\log\a\right)}{z^2+t^2}\, dt,
\end{align}
since $\Gamma(z)=(z-1)!$.
To obtain (\ref{gui}), we simply replace $k$ by $n$, $z-1$ by $k$, $\a$ by $x$ and  $\b$ by $1/x$ in the first equality of (\ref{guigen}).
\section{The second modular relation involving Hurwitz zeta functions}
\begin{theorem}
Let $0<$ Re $z<2$. Define $\varphi(z,x)$ as
\begin{equation}
\varphi(z,x)=\zeta(z,x)-\frac{1}{2}x^{-z}+\frac{x^{1-z}}{1-z},
\end{equation}
where $\zeta(z,x)$ denotes the Hurwitz zeta function. Then if $\alpha$ and $\beta$ are any positive numbers such that $\alpha\beta=1$,
\begin{align}\label{maineq2}
&\a^{\frac{z}{2}}\left(\sum_{n=1}^{\infty}\varphi(z,n\a)-\frac{\zeta(z)}{2\alpha^{z}}-\frac{\zeta(z-1)}{\alpha (z-1)}\right)=\b^{\frac{z}{2}}\left(\sum_{n=1}^{\infty}\varphi(z,n\b)-\frac{\zeta(z)}{2\beta^{z}}-\frac{\zeta(z-1)}{\beta (z-1)}\right)\nonumber\\
&=\frac{8(4\pi)^{\frac{z-4}{2}}}{\Gamma(z)}\int_{0}^{\infty}\Gamma\left(\frac{z-2+it}{4}\right)\Gamma\left(\frac{z-2-it}{4}\right)
\Xi\left(\frac{t+i(z-1)}{2}\right)\Xi\left(\frac{t-i(z-1)}{2}\right)\frac{\cos\left( \tf{1}{2}t\log\a\right)}{z^2+t^2}\, dt,
\end{align}
where $\Xi(t)$ is defined in (\ref{xif}).
\end{theorem}
\begin{proof} 
The asymptotic expansion of $\zeta(z,x)$ \cite[p.~25]{mos} for large $|x|$ and $|\arg$ $x|<\pi$ is given by
\begin{equation}\label{hurae}
\zeta(z,x)=\frac{1}{\Gamma(z)}\left(x^{1-z}\Gamma(z-1)+\frac{1}{2}\Gamma(z)x^{-z}+\sum_{n=1}^{m-1}\frac{B_{2n}}{(2n)!}\Gamma(z+2n-1)x^{-2n-z+1}\right)+O(x^{-2m-z-1}).
\end{equation}
Hence for $0<$ Re $z<2$, the series $\displaystyle\sum_{n=1}^{\infty}\varphi(z,n\a)$ as well as $\displaystyle\sum_{n=1}^{\infty}\varphi(z,n\b)$ are analytic functions. The proof of (\ref{maineq2}) that we give is along the similar lines as the first proof in \cite{bcbad} except that here we have to take care of two parameters $z$ and $\a$ instead of just $\a$ in \cite{bcbad}. We first prove the result for $1<$ Re $z<2$ and later extend it to $0<$ Re $z<2$ using analytic continuation. 

Replacing $s$ by $z-1$ in (\ref{eq20}), we find that for $0<$ Re $z<2$,
\begin{align}\label{eq20modd}
&\int_{0}^{\infty}\Gamma\left(\frac{z-2+it}{4}\right)\Gamma\left(\frac{z-2-it}{4}\right)
\Xi\left(\frac{t+i(z-1)}{2}\right)\Xi\left(\frac{t-i(z-1)}{2}\right)\frac{\cos nt}{z^2+t^2}\, dt\nonumber\\
&=\frac{1}{8}(4\pi)^{-\frac{(z-4)}{2}}\int_{0}^{\infty}x^{z-1}\left(\frac{1}{e^{xe^{n}}-1}-\frac{1}{xe^{n}}\right)\left(\frac{1}{e^{xe^{-n}}-1}-\frac{1}{xe^{-n}}\right)\, dx.
\end{align}
Multiplying both sides of (\ref{eq20modd}) by $\displaystyle8(4\pi)^{\frac{(z-4)}{2}}$ and then letting $n=\displaystyle\tf{1}{2}\log\a$, we see that
\begin{align}\label{eq20modd2}
&8(4\pi)^{\frac{(z-4)}{2}}\int_{0}^{\infty}\Gamma\left(\frac{z-2+it}{4}\right)\Gamma\left(\frac{z-2-it}{4}\right)
\Xi\left(\frac{t+i(z-1)}{2}\right)\Xi\left(\frac{t-i(z-1)}{2}\right)\frac{\cos \left(\tf{1}{2}t\log\a\right)}{z^2+t^2}\, dt\nonumber\\
&=\int_{0}^{\infty}x^{z-1}\left(\frac{1}{e^{x\sqrt{\a}}-1}-\frac{1}{x\sqrt{\a}}\right)\left(\frac{1}{e^{x/\sqrt{\a}}-1}-\frac{1}{x/\sqrt{\a}}\right)\, dx.
\end{align}
Making a change of variable $x=2\pi t/\sqrt{\a}$ in the integral on the right-hand side of (\ref{eq20modd}), we see that
\begin{align}\label{eq20modd2n}
&8(4\pi)^{\frac{(z-4)}{2}}\int_{0}^{\infty}\Gamma\left(\frac{z-2+it}{4}\right)\Gamma\left(\frac{z-2-it}{4}\right)
\Xi\left(\frac{t+i(z-1)}{2}\right)\Xi\left(\frac{t-i(z-1)}{2}\right)\frac{\cos \left(\tf{1}{2}t\log\a\right)}{z^2+t^2}\, dt\nonumber\\
&=\left(\frac{2\pi}{\sqrt{\a}}\right)^{z}\int_{0}^{\infty}t^{z-1}\left(\frac{1}{e^{2\pi t}-1}-\frac{1}{2\pi t}\right)\left(\frac{1}{e^{2\pi t/\a}-1}-\frac{1}{2\pi t/\a}\right)\, dt\nonumber\\
&=\left(\frac{2\pi}{\sqrt{\a}}\right)^{z}\int_{0}^{\infty}t^{z-1}\left(\frac{1}{(e^{2\pi
t/\alpha}-1)(e^{2\pi t}-1)} -\frac{\a}{2\pi t(e^{2\pi
t}-1)}-\frac{1}{2\pi t(e^{2\pi t/\alpha}-1)}
+\frac{\alpha}{4\pi^2 t^{2}}\right)\, dt\nonumber\\
&=I_{1}(z,\a)+I_{2}(z,\a),
\end{align}
where 
\begin{equation}\label{inti1}
I_{1}(z,\a)=\left(\frac{2\pi}{\sqrt{\a}}\right)^{z}\int_{0}^{\infty}t^{z-1}\left(\frac{1}{(e^{2\pi
t/\alpha}-1)(e^{2\pi t}-1)} -\frac{\alpha}{2\pi t(e^{2\pi
t}-1)}+\frac{e^{-t/\a}}{4\pi t}\right)\, dt,
\end{equation}
and
\begin{equation}\label{inti2}
I_{2}(z,\a)=\left(\frac{2\pi}{\sqrt{\a}}\right)^{z}\int_{0}^{\infty}t^{z-1}\left(\frac{-1}{2\pi t(e^{2\pi t/\alpha}-1)}
+\frac{\alpha}{4\pi^2 t^{2}}-\frac{e^{-t/\a}}{4\pi t}\right)\, dt.
\end{equation}

First,
\begin{align}\label{cali1}
I_{1}(z,\a)&=\left(\frac{2\pi}{\sqrt{\a}}\right)^{z}\times\nonumber\\
\qquad &\int_{0}^{\infty}t^{z-1}\left(\frac{1}{(e^{2\pi
t/\alpha}-1)(e^{2\pi t}-1)}-\frac{\alpha}{2\pi t(e^{2\pi t}-1)}+\frac{1}{2(e^{2\pi t}-1)}-\frac{1}{2(e^{2\pi t}-1)}+\frac{e^{-t/\a}}{4\pi t}\right)\, dt\nonumber\\
&=\left(\frac{2\pi}{\sqrt{\a}}\right)^{z}\int_{0}^{\infty}\frac{t^{z-1}}{(e^{2\pi t}-1)}\left(\frac{1}{(e^{2\pi t/\alpha}-1)}-\frac{1}{2\pi t/\a}+\frac{1}{2}\right)\, dt\nonumber\\
&\quad-\frac{(2\pi)^{z-1}}{2\a^{z/2}}\int_{0}^{\infty}t^{z-1}\left(\frac{2\pi}{e^{2\pi t}-1}-\frac{e^{-t/\a}}{t}\right)\, dt\nonumber\\
&=I_{3}(z,\a)+I_{4}(z,\a),
\end{align}
where 
\begin{equation}\label{inti3}
I_{3}(z,\a)=\left(\frac{2\pi}{\sqrt{\a}}\right)^{z}\int_{0}^{\infty}\frac{t^{z-1}}{(e^{2\pi t}-1)}\left(\frac{1}{(e^{2\pi t/\alpha}-1)}-\frac{1}{2\pi t/\a}+\frac{1}{2}\right)\, dt,
\end{equation}
and
\begin{equation}\label{inti4}
I_{4}(z,\a)=-\frac{(2\pi)^{z-1}}{2\a^{z/2}}\int_{0}^{\infty}t^{z-1}\left(\frac{2\pi}{e^{2\pi t}-1}-\frac{e^{-t/\a}}{t}\right)\, dt.
\end{equation}
Now for evaluating $I_{3}(z,\a)$, we make use of the following formula \cite[p.~23]{mos} valid for Re $z>-1$ and Re $a>0$
\begin{equation}\label{melvarphi}
\zeta(z,a)=\frac{a^{-z}}{2}-\frac{a^{1-z}}{1-z}+\frac{1}{\Gamma(z)}\int_{0}^{\infty}e^{-ax}x^{z-1}\left(\frac{1}{e^{x}-1}-\frac{1}{x}+\frac{1}{2}\right)\, dx.
\end{equation}
Since $t>0$, expanding $\displaystyle\frac{1}{(e^{2\pi t}-1)}$ in terms of its geometric series, interchanging the summation and integration because of absolute convergence and then using (\ref{melvarphi}), we find that
\begin{align}\label{cali3}
I_{3}(z,\a)&=\left(\frac{2\pi}{\sqrt{\a}}\right)^{z}\int_{0}^{\infty}\frac{t^{z-1}e^{-2\pi t}}{1-e^{-2\pi t}}\left(\frac{1}{(e^{2\pi t/\alpha}-1)}-\frac{1}{2\pi t/\a}+\frac{1}{2}\right)\, dt\nonumber\\
&=\left(\frac{2\pi}{\sqrt{\a}}\right)^{z}\sum_{n=1}^{\infty}\int_{0}^{\infty}t^{z-1}e^{-2\pi nt}\left(\frac{1}{(e^{2\pi t/\alpha}-1)}-\frac{1}{2\pi t/\a}+\frac{1}{2}\right)\, dt\nonumber\\
&=\a^{z/2}\sum_{n=1}^{\infty}\Gamma(z)\left(\zeta(z,n\a)-\frac{(n\a)^{-z}}{2}+\frac{(n\a)^{1-z}}{1-z}\right)\nonumber\\
&=\a^{z/2}\Gamma(z)\sum_{n=1}^{\infty}\varphi(z,n\a),
\end{align}
where in the antepenultimate line, we have made a change of variable $x=2\pi t/\a$. Next we evaluate $I_{4}(z,\a)$. Since Re $z>1$, using (\ref{gamzeta}) and the integral representation for $\Gamma(z-1)$, we find that
\begin{align}\label{cali4}
I_{4}(z,\a)&=-\frac{(2\pi)^{z-1}}{2\a^{z/2}}\left(2\pi\int_{0}^{\infty}\frac{t^{z-1}}{e^{2\pi t}-1}\, dt-\int_{0}^{\infty}t^{z-2}e^{-t/\a}\, dt\right)\nonumber\\
&=-\frac{(2\pi)^{z-1}}{2\a^{z/2}}\left(\frac{2\pi}{(2\pi)^{z}}\Gamma(z)\zeta(z)-\a^{z-1}\Gamma(z-1)\right)\nonumber\\
&=-\frac{\a^{-\tf{z}{2}}}{2}\Gamma(z)\zeta(z)+\frac{\a^{\tf{z}{2}-1}}{2}(2\pi)^{z-1}\Gamma(z-1).
\end{align}
Hence from (\ref{cali1}), (\ref{cali3}) and (\ref{cali4}), we conclude that
\begin{equation}\label{finali1}
I_{1}(z,\a)=\a^{z/2}\Gamma(z)\sum_{n=1}^{\infty}\varphi(z,n\a)-\frac{\a^{-\tf{z}{2}}}{2}\Gamma(z)\zeta(z)+\frac{\a^{\tf{z}{2}-1}}{2}(2\pi)^{z-1}\Gamma(z-1)
\end{equation}

It remains to evaluate $I_{2}(z,\a)$. Now from \cite[p.~23, eqn. 2.7.1]{titch} for $0<$ Re $z<1$, we have
\begin{equation}\label{gammazeta}
\Gamma(z)\zeta(z)=\int_{0}^{\infty}t^{z-1}\left(\frac{1}{e^{t}-1}-\frac{1}{t}\right)\, dt
\end{equation}
Thus employing a change of variable $u=2\pi t/\a$ in (\ref{inti2}) and then using (\ref{gammazeta}) and the integral representation for $\Gamma(z-1)$, we see that
\begin{align}\label{finali2}
I_{2}(z,\a)&=-\a^{\tf{z}{2}-1}\int_{0}^{\infty}u^{z-1}\left(\frac{1}{u(e^{u}-1)}-\frac{1}{u^2}+\frac{e^{-u/2\pi}}{2u}\right)\, du\nonumber\\
&=-\a^{\tf{z}{2}-1}\left(\int_{0}^{\infty}u^{z-2}\left(\frac{1}{e^{u}-1}-\frac{1}{u}\right)\, du+\frac{1}{2}\int_{0}^{\infty}e^{-u/2\pi}u^{z-2}\, du\right)\nonumber\\
&=-\a^{\tf{z}{2}-1}\left(\Gamma(z-1)\zeta(z-1)+\frac{1}{2}(2\pi)^{z-1}\Gamma(z-1)\right).
\end{align}
Finally from (\ref{eq20modd2n}), (\ref{finali1}) and (\ref{finali2}), we see that after simplification
\begin{align}\label{eq20modd2n1}
&\frac{8(4\pi)^{\frac{(z-4)}{2}}}{\Gamma(z)}\int_{0}^{\infty}\Gamma\left(\frac{z-2+it}{4}\right)\Gamma\left(\frac{z-2-it}{4}\right)
\Xi\left(\frac{t+i(z-1)}{2}\right)\Xi\left(\frac{t-i(z-1)}{2}\right)\frac{\cos \left(\tf{1}{2}t\log\a\right)}{z^2+t^2}\, dt\nonumber\\
&=\a^{\frac{z}{2}}\left(\sum_{n=1}^{\infty}\varphi(z,n\a)-\frac{\zeta(z)}{2\alpha^{z}}-\frac{\zeta(z-1)}{\alpha (z-1)}\right).
\end{align}
Now replacing $\a$ by $\b$ in (\ref{eq20modd2n1}), we find that
\begin{align}\label{thmover}
&\b^{\frac{z}{2}}\left(\sum_{n=1}^{\infty}\varphi(z,n\b)-\frac{\zeta(z)}{2\beta^{z}}-\frac{\zeta(z-1)}{\beta (z-1)}\right)\nonumber\\
&=\frac{8(4\pi)^{\frac{(z-4)}{2}}}{\Gamma(z)}\int_{0}^{\infty}\Gamma\left(\frac{z-2+it}{4}\right)\Gamma\left(\frac{z-2-it}{4}\right)
\Xi\left(\frac{t+i(z-1)}{2}\right)\Xi\left(\frac{t-i(z-1)}{2}\right)\frac{\cos\left( \tf{1}{2}t\log\b\right)}{z^2+t^2}\, dt\nonumber\\
&=\frac{8(4\pi)^{\frac{(z-4)}{2}}}{\Gamma(z)}\int_{0}^{\infty}\Gamma\left(\frac{z-2+it}{4}\right)\Gamma\left(\frac{z-2-it}{4}\right)
\Xi\left(\frac{t+i(z-1)}{2}\right)\Xi\left(\frac{t-i(z-1)}{2}\right)\frac{\cos\left( \tf{1}{2}t\log\tf{1}{\a}\right)}{z^2+t^2}\, dt\nonumber\\
&=\frac{8(4\pi)^{\frac{(z-4)}{2}}}{\Gamma(z)}\int_{0}^{\infty}\Gamma\left(\frac{z-2+it}{4}\right)\Gamma\left(\frac{z-2-it}{4}\right)
\Xi\left(\frac{t+i(z-1)}{2}\right)\Xi\left(\frac{t-i(z-1)}{2}\right)\frac{\cos\left( \tf{1}{2}t\log\a\right)}{z^2+t^2}\, dt.\nonumber\\
\end{align}
Thus (\ref{eq20modd2n1}) and (\ref{thmover}) imply (\ref{maineq2}) for $1<$ Re $z<2$. Now observe from (\ref{eq20modd}) that the right-hand side of (\ref{maineq2}) is analytic for $0<$ Re $z<2$. Also on the left-hand side of (\ref{maineq2}), $z=1$ is a removable singularity because the residue of $\zeta(z)$ at $z=1$ is equal to $1$ and because $\zeta(0)=-\frac{1}{2}$. Hence by analytic continuation, (\ref{maineq2}) holds for $0<$ Re $z<2$.
\end{proof}

As a limiting case of (\ref{maineq2}), we obtain Ramanujan's transformation formula, i.e., equation (\ref{w1.26}).
\begin{corollary}
If
\begin{equation}\label{psi2def}
\phi(x):=\psi(x)+\df{1}{2x}-\log x,
\end{equation}
and $\a$ and $\b$ are positive numbers such that $\a\b=1$, then
\begin{multline}\label{ramcor}
\sqrt{\a}\left\{\df{\gamma-\log(2\pi\a)}{2\a}+\sum_{n=1}^{\infty}\phi(n\a)\right\}
=\sqrt{\b}\left\{\df{\gamma-\log(2\pi\b)}{2\b}+\sum_{n=1}^{\infty}\phi(n\b)\right\}\\
=-\df{1}{\pi^{3/2}}\int_0^{\infty}\left|\Xi\left(\df{1}{2}t\right)\Gamma\left(\df{-1+it}{4}\right)\right|^2
\df{\cos\left(\tf{1}{2}t\log\a\right)}{1+t^2}\, dt.
\end{multline}
\end{corollary}
\begin{proof}
Let $z\to 1$ in (\ref{maineq2}). Then using Lebesgue's dominated convergence theorem, we observe that
\begin{align}\label{maineq2lim}
&\lim_{z\to 1}\a^{\frac{z}{2}}\left(\sum_{n=1}^{\infty}\varphi(z,n\a)-\frac{\zeta(z)}{2\alpha^{z}}-\frac{\zeta(z-1)}{\alpha (z-1)}\right)=\lim_{z\to 1}\b^{\frac{z}{2}}\left(\sum_{n=1}^{\infty}\varphi(z,n\b)-\frac{\zeta(z)}{2\beta^{z}}-\frac{\zeta(z-1)}{\beta (z-1)}\right)\nonumber\\
&=\df{1}{\pi^{3/2}}\int_0^{\infty}\left|\Xi\left(\df{1}{2}t\right)\Gamma\left(\df{-1+it}{4}\right)\right|^2
\df{\cos\left(\tf{1}{2}t\log\a\right)}{1+t^2}\, dt
\end{align}
Now since $\sum_{n=1}^{\infty}\varphi(z,n\alpha)$ and $\sum_{n=1}^{\infty}\varphi(z,n\beta)$ converge absolutely and uniformly in a neighborhood of $z=1$, which can be seen from (\ref{hurae}), we observe that
\begin{align}\label{maineq2limsim}
&\lim_{z\to 1}\a^{\frac{z}{2}}\left(\sum_{n=1}^{\infty}\varphi(z,n\a)-\frac{\zeta(z)}{2\alpha^{z}}-\frac{\zeta(z-1)}{\alpha (z-1)}\right)\nonumber\\
&=\left(\lim_{z\to 1}\alpha^{\frac{z}{2}}\right)\cdot\left(\sum_{n=1}^{\infty}\left(\lim_{z\to 1}\varphi(z,n\a)\right)-\lim_{z\to 1}\left(\frac{\zeta(z)}{2\alpha^{z}}+\frac{\zeta(z-1)}{\alpha (z-1)}\right)\right)\nonumber\\
&=\left(\lim_{z\to 1}\alpha^{\frac{z}{2}}\right)\cdot\left(\sum_{n=1}^{\infty}\left(\lim_{z\to 1}\varphi(z,n\a)\right)-\lim_{z\to 1}\left(\frac{\zeta(z)}{2\alpha^{z}}-\frac{1}{2\alpha^{z}(z-1)}\right)-\lim_{z\to 1}\left(\frac{1}{2\alpha^{z}(z-1)}+\frac{\zeta(z-1)}{\alpha (z-1)}\right)\right)\nonumber\\
\end{align}
But it is known \cite[p.~16]{titch} that
\begin{equation}\label{zetaprop}
\lim_{s\to 1}\left(\zeta(s)-\frac{1}{s-1}\right)=\gamma.
\end{equation}
Hence
\begin{equation}\label{seclim}
\lim_{z\to 1}\left(\frac{\zeta(z)}{2\alpha^{z}}-\frac{1}{2\alpha^{z}(z-1)}\right)=\frac{\gamma}{2\alpha}.
\end{equation}
Next using L'Hopital's rule, we see that
{\allowdisplaybreaks\begin{align}\label{thlim}
\lim_{z\to 1}\left(\frac{1}{2\alpha^{z}(z-1)}+\frac{\zeta(z-1)}{\alpha (z-1)}\right)&=\lim_{z\to 1}\frac{1}{2\alpha^{z}}\cdot
\lim_{z\to 1}\frac{1+2\alpha^{z-1}\zeta(z-1)}{z-1}\nonumber\\
&=\lim_{z\to 1}\frac{1}{2\alpha^{z}}\cdot\lim_{z\to 1}\left(2\alpha^{z-1}\zeta^{'}(z-1)+2\zeta(z-1)\alpha^{z-1}\log\alpha\right)\nonumber\\
&=\frac{1}{2\alpha}\left(2\zeta^{'}(0)+2\zeta(0)\log\alpha\right)\nonumber\\
&=-\frac{\log(2\pi\alpha)}{2\alpha},
\end{align}}
since $\zeta(0)=-\frac{1}{2}$ and $\zeta^{'}(0)=-\frac{1}{2}\log(2\pi)$ \cite[pp.~19-20, eqns. 2.4.3, 2.4.5]{titch}.
Now noting that \cite[p.~23]{mos}
\begin{equation}\label{blerch}
\lim_{z\to 1}\left(\zeta(z,a)-\frac{1}{z-1}\right)=-\psi(a),
\end{equation}
 and using L'Hopital's rule again, we observe that
{\allowdisplaybreaks\begin{align}\label{firlim}
&\lim_{z\to 1}\varphi(n\a)\nonumber\\
&=\lim_{z\to 1}\left(\zeta(z,n\alpha)-\frac{1}{2}(n\alpha)^{-z}+\frac{(n\alpha)^{1-z}}{1-z}\right)\nonumber\\
&=\lim_{z\to 1}\left[\left(\zeta(z,n\alpha)-\frac{1}{z-1}\right)-\frac{(n\alpha)^{-z}}{2}+\frac{(n\alpha)^{1-z}-1}{1-z}\right]\nonumber\\
&=-\psi(n\alpha)-\frac{1}{2n\alpha}+\lim_{z\to 1}\frac{-(n\alpha)^{1-z}\log n\alpha}{-1}\nonumber\\
&=-\psi(n\alpha)-\frac{1}{2n\alpha}+\log (n\alpha)\nonumber\\
&=-\phi(n\alpha),
\end{align}}
where $\phi(x)$ is as defined in (\ref{psi2def}).
Hence from (\ref{maineq2limsim}), (\ref{seclim}), (\ref{thlim}) and (\ref{firlim}), we find that
\begin{equation}\label{fin1side}
\lim_{z\to 1}\a^{\frac{z}{2}}\left(\sum_{n=1}^{\infty}\varphi(z,n\a)-\frac{\zeta(z)}{2\alpha^{z}}-\frac{\zeta(z-1)}{\alpha (z-1)}\right)=-\sqrt{\alpha}\left(\frac{\gamma-\log(2\pi\alpha)}{2\alpha}+\sum_{n=1}^{\infty}\phi(n\alpha)\right),
\end{equation}
Thus from (\ref{maineq2lim}), (\ref{fin1side}) and (\ref{fin1side}) with $\alpha$ replaced by $\beta$, we obtain 
\begin{multline}\label{ramcor2}
-\sqrt{\a}\left\{\df{\gamma-\log(2\pi\a)}{2\a}+\sum_{n=1}^{\infty}\phi(z,n\a)\right\}
=-\sqrt{\b}\left\{\df{\gamma-\log(2\pi\b)}{2\b}+\sum_{n=1}^{\infty}\phi(z,n\b)\right\}\\
=\df{1}{\pi^{3/2}}\int_0^{\infty}\left|\Xi\left(\df{1}{2}t\right)\Gamma\left(\df{-1+it}{4}\right)\right|^2
\df{\cos\left(\tf{1}{2}t\log\a\right)}{1+t^2}\, dt.
\end{multline}
Multiplying (\ref{ramcor2}) throughout by $-1$, we arrive at (\ref{ramcor}). 
\end{proof}
\section{Proof of (\ref{eq19cor})}
Here we show that the second term on the right-hand side of equation (19) in \cite{riemann}, namely, $\displaystyle-\frac{1}{4}(4\pi)^{\frac{(s-3)}{2}}\Gamma(s)\zeta(s)\cosh n(1-s)$  is not present and thus indeed that (\ref{eq19cor}) is actually the correct version of equation (19) in \cite{riemann}. Since the exposition in Sections $4$ and $5$ of \cite{riemann} is quite terse, we will derive (\ref{eq19cor}) giving all the details. We will collect and prove wherever necessary, several ingredients required for the proof along the way.

First, equation (15) in \cite{riemann} states that for Re $s>-1$ and $\alpha\beta=4\pi^2$, we have
\begin{align}\label{modrel}
G(s ;\alpha)&:=\frac{\zeta(1-s)}{4\cos(\pi s/2)}\alpha^{(s-1)/2}+\frac{\zeta(-s)}{8\sin(\pi s/2)}\alpha^{(s+1)/2}+\alpha^{(s+1)/2}\int_{0}^{\infty}\int_{0}^{\infty}\frac{x^s \sin\left(\alpha xy\right)}{(e^{2\pi x}-1)(e^{2\pi y}-1)}\, dx dy\nonumber\\
&=\frac{\zeta(1-s)}{4\cos(\pi s/2)}\beta^{(s-1)/2}+\frac{\zeta(-s)}{8\sin(\pi s/2)}\beta^{(s+1)/2}+\beta^{(s+1)/2}\int_{0}^{\infty}\int_{0}^{\infty}\frac{x^s \sin\left(\beta xy\right)}{(e^{2\pi x}-1)(e^{2\pi y}-1)}\, dx dy.
\end{align}
This relation can be proved by obtaining integral representations for $\displaystyle\frac{\zeta(1-s)}{4\cos(\pi s/2)}$ and $\displaystyle\frac{\zeta(-s)}{8\sin(\pi s/2)}$ and by using the identity (see \cite[p.~253]{riemann})
\begin{equation}\label{sine}
\int_{0}^{\infty}\frac{\sin(\alpha xy)}{e^{2\pi y}-1}\, dy=\frac{1}{2}\left(\frac{1}{e^{\alpha x}-1}-\frac{1}{\alpha x}+\frac{1}{2}\right).
\end{equation}
Since we are concerned with the case Re $s>1$, we prove (\ref{modrel}) for Re $s>1$ only. Other cases can be similarly proved.

Now the functional equation for $\zeta(s)$ in its non-symmetric form \cite[p.~13, eqn. 2.1.1]{titch} states that
\begin{equation}\label{zetafe}
\zeta(s)=2^s\pi^{s-1}\Gamma(1-s)\zeta(1-s)\sin\left(\tf{1}{2}\pi s\right).
\end{equation}
Using this and (\ref{gamzeta}), one can easily show that for Re $s>1$,
\begin{equation}\label{int1}
\frac{\zeta(1-s)}{4\cos(\tfrac{1}{2}\pi s)}=\frac{1}{2}\int_{0}^{\infty}\frac{x^{s-1} dx}{e^{2\pi x}-1},
\end{equation}
and
\begin{equation}\label{int2}
\frac{\zeta(-s)}{8\sin(\tfrac{1}{2}\pi s)}=\frac{-1}{4}\int_{0}^{\infty}\frac{x^{s} dx}{e^{2\pi x}-1}.
\end{equation}
Hence using (\ref{sine}), (\ref{int1}) and (\ref{int2}) in (\ref{modrel}), we see that
\begin{equation}\label{modleft}
G(s;\alpha)=\frac{\alpha^{(s+1)/2}}{2}\int_{0}^{\infty}\frac{x^s}{(e^{2\pi x}-1)(e^{\alpha x}-1)}\, dx.
\end{equation}
So (\ref{modrel}) will be proved for Re $s>1$ if we can show that
\begin{equation}\label{modlr}
\frac{\alpha^{(s+1)/2}}{2}\int_{0}^{\infty}\frac{x^s}{(e^{2\pi x}-1)(e^{\alpha x}-1)}\, dx=\frac{\beta^{(s+1)/2}}{2}\int_{0}^{\infty}\frac{x^s}{(e^{2\pi x}-1)(e^{\beta x}-1)}\, dx.
\end{equation}
But this is easily seen by making the substitution $x=\displaystyle\frac{2\pi y}{\alpha}$ on the left-hand side of (\ref{modlr}) and using the fact that
$\alpha\beta=4\pi^2$. Thus (\ref{modrel}) is proved for Re $s>1$.

Next, equation (17) in \cite{riemann} states that when $\alpha\beta=4\pi^2$ and Re $s>-1$, we have
\begin{align}\label{17ram}
&\frac{\zeta(1-s)}{4\cos(\pi s/2)}\frac{s-1}{(s-1)^2+t^2}\left(\alpha^{(s-1)/2}+\beta^{(s-1)/2}\right)+\frac{\zeta(-s)}{8\sin(\pi s/2)}\frac{s+1}{(s+1)^2+t^2}\left(\alpha^{(s+1)/2}+\beta^{(s+1)/2}\right)\nonumber\\
&+\alpha^{(s+1)/2}\int_{0}^{\infty}\int_{0}^{\infty}\left(\frac{\alpha xy}{1!}\frac{s+3}{(s+3)^2+t^2}-\frac{(\alpha xy)^3}{3!}\frac{s+7}{(s+7)^2+t^2}+\cdots\right)\frac{x^s dx dy}{(e^{2\pi x}-1)(e^{2\pi y}-1)}\nonumber\\
&+\beta^{(s+1)/2}\int_{0}^{\infty}\int_{0}^{\infty}\left(\frac{\beta xy}{1!}\frac{s+3}{(s+3)^2+t^2}-\frac{(\beta xy)^3}{3!}\frac{s+7}{(s+7)^2+t^2}+\cdots\right)\frac{x^s dx dy}{(e^{2\pi x}-1)(e^{2\pi y}-1)}\nonumber\\
&=\frac{2^{(s-3)/2}}{\pi}\frac{\Gamma\left(\tf{1}{4}(s-1+it)\right)\Gamma\left(\tf{1}{4}(s-1-it)\right)}{(s+1)^2+t^2}\Xi\left(\frac{t+is}{2}
\right)\Xi\left(\frac{t-is}{2}\right)\cos\left(\frac{1}{4}t\log\frac{\alpha}{\beta}\right).
\end{align}

Letting $\alpha=\beta=2\pi$ in (\ref{17ram}) and simplifying, we find that
\begin{align}\label{17ram2}
&\frac{\zeta(1-s)}{\pi\cos(\pi s/2)}\frac{s-1}{(s-1)^2+t^2}+\frac{\zeta(-s)}{\sin(\pi s/2)}\frac{s+1}{(s+1)^2+t^2}\nonumber\\
&+8\int_{0}^{\infty}\int_{0}^{\infty}\left(\frac{2\pi xy}{1!}\frac{s+3}{(s+3)^2+t^2}-\frac{(2\pi xy)^3}{3!}\frac{s+7}{(s+7)^2+t^2}+\cdots\right)\frac{x^s dx dy}{(e^{2\pi x}-1)(e^{2\pi y}-1)}\nonumber\\
&=\frac{1}{\pi^{(s+3)/2}}\frac{\Gamma\left(\tf{1}{4}(s-1+it)\right)\Gamma\left(\tf{1}{4}(s-1-it)\right)}{(s+1)^2+t^2}\Xi\left(\frac{t+is}{2}
\right)\Xi\left(\frac{t-is}{2}\right).
\end{align}
Now we know that for Re $a>0$,
\begin{equation}\label{use}
\int_{0}^{\infty}e^{-au}\cos bu\, du=\frac{a}{a^2+b^2}.
\end{equation}

Using (\ref{use}), we wish to replace the fractions of the form $\displaystyle\frac{s+j}{(s+j)^2+t^2}$ in (\ref{17ram2}) by integrals. Since Re $s>1$, for $j\geq -1$, we have
\begin{equation}\label{useint}
\frac{s+j}{(s+j)^2+t^2}=\int_{0}^{\infty}e^{-(s+j)u}\cos tu\, du.
\end{equation}

Hence using (\ref{useint}) in (\ref{17ram2}), inverting the order of integration because of absolute convergence and simplifying, we see that
\begin{align}\label{17ram3}
&\frac{1}{\pi^{(s+3)/2}}\frac{\Gamma\left(\tf{1}{4}(s-1+it)\right)\Gamma\left(\tf{1}{4}(s-1-it)\right)}{(s+1)^2+t^2}\Xi\left(\frac{t+is}{2}
\right)\Xi\left(\frac{t-is}{2}\right)\nonumber\\
&=\frac{\zeta(1-s)}{\pi\cos(\pi s/2)}\int_{0}^{\infty}e^{-(s-1)u}\cos tu\, du+\frac{\zeta(-s)}{\sin(\pi s/2)}\int_{0}^{\infty}e^{-(s+1)u}\cos tu\, du\nonumber\\
&\quad+8\int_{0}^{\infty}\int_{0}^{\infty}\left(\frac{2\pi xy}{1!}\int_{0}^{\infty}e^{-(s+3)u}\cos tu\, du-\frac{(2\pi xy)^3}{3!}\int_{0}^{\infty}e^{-(s+7)u}\cos tu\, du+\cdots\right)\nonumber\\
&\quad\quad\quad\times\frac{x^s dx dy}{(e^{2\pi x}-1)(e^{2\pi y}-1)}\nonumber\\
&=\int_{0}^{\infty}\bigg[\frac{\zeta(1-s)}{\pi\cos(\pi s/2)}e^{-(s-1)u}+\frac{\zeta(-s)}{\sin(\pi s/2)}e^{-(s+1)u}\nonumber\\
&\quad+8e^{-(s+1)u}\int_{0}^{\infty}\int_{0}^{\infty}\left(\frac{2\pi xye^{-2u}}{1!}-\frac{(2\pi xye^{-2u})^3}{3!}+\cdots\right)\frac{x^s dx dy}{(e^{2\pi x}-1)(e^{2\pi y}-1)}\bigg]\cos tu\, du\nonumber\\
&=\int_{0}^{\infty}\bigg[\frac{\zeta(1-s)}{\pi\cos(\pi s/2)}e^{-(s-1)u}+\frac{\zeta(-s)}{\sin(\pi s/2)}e^{-(s+1)u}+8e^{-(s+1)u}\int_{0}^{\infty}\int_{0}^{\infty}\frac{x^s \sin\left(2\pi xye^{-2u}\right)}{(e^{2\pi x}-1)(e^{2\pi y}-1)}\, dx dy\bigg]\nonumber\\
&\quad\quad\quad\times\cos tu\, du.
\end{align}

Now let
\begin{equation}\label{deff}
f(u):=\frac{\zeta(1-s)}{\pi\cos(\pi s/2)}e^{-(s-1)u}+\frac{\zeta(-s)}{\sin(\pi s/2)}e^{-(s+1)u}+8e^{-(s+1)u}\int_{0}^{\infty}\int_{0}^{\infty}\frac{x^s \sin\left(2\pi xye^{-2u}\right)}{(e^{2\pi x}-1)(e^{2\pi y}-1)}\, dx dy,
\end{equation}
and
\begin{equation}\label{deff2}
\widehat{f}(t):=\frac{1}{\pi^{(s+3)/2}}\frac{\Gamma\left(\tf{1}{4}(s-1+it)\right)\Gamma\left(\tf{1}{4}(s-1-it)\right)}{(s+1)^2+t^2}\Xi\left(\frac{t+is}{2}
\right)\Xi\left(\frac{t-is}{2}\right).
\end{equation}
Then from (\ref{17ram3}), (\ref{deff}) and (\ref{deff2}), we have
\begin{equation}\label{relff2}
\widehat{f}(t)=\int_{0}^{\infty}f(u)\cos tu\, du.
\end{equation}
Now we show that $f$ is an even function of $u$.\\

If we let $\alpha=2\pi e^{-2u}$ and $\beta=2\pi e^{2u}$ in (\ref{modrel}), upon simplification, we find that
\begin{align}\label{feven}
&\frac{\zeta(1-s)}{\pi\cos(\pi s/2)}e^{-(s-1)u}+\frac{\zeta(-s)}{\sin(\pi s/2)}e^{-(s+1)u}+8e^{-(s+1)u}\int_{0}^{\infty}\int_{0}^{\infty}\frac{x^s \sin\left(2\pi xye^{-2u}\right)}{(e^{2\pi x}-1)(e^{2\pi y}-1)}\, dx dy\nonumber\\
&=\frac{\zeta(1-s)}{\pi\cos(\pi s/2)}e^{(s-1)u}+\frac{\zeta(-s)}{\sin(\pi s/2)}e^{(s+1)u}+8e^{(s+1)u}\int_{0}^{\infty}\int_{0}^{\infty}\frac{x^s \sin\left(2\pi xye^{2u}\right)}{(e^{2\pi x}-1)(e^{2\pi y}-1)}\, dx dy.
\end{align}
This proves that $f$ is an even function of $u$. Also using the fact that $\Xi(-t)=\Xi(t)$, we readily observe that $\widehat{f}$ is also an even function of $t$.

Then from Fourier's integral theorem and (\ref{relff2}), we deduce that for $n$ real,
\begin{equation}\label{relfou}
f(n)=\frac{2}{\pi}\int_{0}^{\infty}\widehat{f}(t)\cos nt\, dt.
\end{equation}
Now define
\begin{equation}
F(n,s):=\displaystyle\frac{\pi^{\frac{s+5}{2}}}{2}f(n).
\end{equation}
Then from (\ref{deff}), (\ref{deff2}) and (\ref{relfou}), we find that
\begin{align}\label{mainrel}
F(n,s)&=\int_{0}^{\infty}\frac{\Gamma\left(\tf{1}{4}(s-1+it)\right)\Gamma\left(\tf{1}{4}(s-1-it)\right)}{(s+1)^2+t^2}\Xi\left(\frac{t+is}{2}
\right)\Xi\left(\frac{t-is}{2}\right)\cos nt\, dt\nonumber\\
&=\frac{\pi^{\frac{s+5}{2}}}{2}\bigg(\frac{\zeta(1-s)}{\pi\cos(\pi s/2)}e^{-(s-1)n}+\frac{\zeta(-s)}{\sin(\pi s/2)}e^{-(s+1)n}+8e^{-(s+1)n}\int_{0}^{\infty}\int_{0}^{\infty}\frac{t^s \sin\left(2\pi tye^{-2n}\right)}{(e^{2\pi t}-1)(e^{2\pi y}-1)}\, dt dy\bigg).
\end{align}
Substituting (\ref{sine}), (\ref{int1}) and (\ref{int2}) in (\ref{mainrel}), we find that
\begin{align}\label{mainrel3}
F(n,s)&=\frac{\pi^{\frac{s+5}{2}}}{2}\bigg[\frac{2e^{-(s-1)n}}{\pi}\int_{0}^{\infty}\frac{t^{s-1} dt}{e^{2\pi t}-1}-2e^{-(s+1)n}
\int_{0}^{\infty}\frac{t^{s} dt}{e^{2\pi t}-1}\nonumber\\
&\quad+4e^{-(s+1)n}\int_{0}^{\infty}\frac{t^s dt}{e^{2\pi t}-1}\left(\frac{1}{e^{2\pi te^{-2n}}-1}-\frac{1}{2\pi te^{-2n}}+\frac{1}{2}\right)\bigg].
\end{align}
Letting $t=\displaystyle\frac{e^n x}{2\pi}$ in (\ref{mainrel3}), we find that
\begin{align}\label{mainrel4}
F(n,s)&=\frac{\pi^{\frac{s+5}{2}}}{2}\bigg[\frac{e^n}{2^{s-1}\pi^{s+1}}\int_{0}^{\infty}\frac{x^{s-1}\, dx}{(e^{xe^n}-1)}-\frac{1}{2^{s}\pi^{s+1}}
\int_{0}^{\infty}\frac{x^{s}\, dx}{(e^{xe^n}-1)}\nonumber\\
&\quad+\frac{1}{2^{s-1}\pi^{s+1}}\int_{0}^{\infty}\frac{x^s dx}{(e^{xe^n}-1)}\left(\frac{1}{e^{xe^{-n}}-1}-\frac{1}{xe^{-n}}+\frac{1}{2}\right)\bigg]\nonumber\\
&=\frac{1}{8}(4\pi)^{-\frac{(s-3)}{2}}\bigg[e^n\int_{0}^{\infty}\frac{x^{s-1}\, dx}{(e^{xe^n}-1)}-\frac{1}{2}
\int_{0}^{\infty}\frac{x^{s}\, dx}{(e^{xe^n}-1)}\nonumber\\
&\quad+\int_{0}^{\infty}\frac{x^s dx}{(e^{xe^n}-1)}\left(\frac{1}{e^{xe^{-n}}-1}-\frac{1}{xe^{-n}}+\frac{1}{2}\right)\bigg]\nonumber\\
&=\frac{1}{8}(4\pi)^{-\frac{(s-3)}{2}}\int_{0}^{\infty}\frac{x^s\, dx}{(e^{xe^n}-1)(e^{xe^{-n}}-1)}.
\end{align}
Finally, we obtain (\ref{eq19cor}) from (\ref{mainrel}) and (\ref{mainrel4}).\\

\textbf{Remark.} 
Equation (\ref{eq20}), for $-1<$ Re $s<1$, is derived in a very similar manner as above, except that since Re $s<1$, the first expression on the left-hand side of (\ref{17ram2}) is written as $\displaystyle-\frac{\zeta(1-s)}{\pi\cos(\pi s/2)}\frac{1-s}{(1-s)^2+t^2}$ and then we use the following equation,
\begin{equation}\label{change}
\frac{1-s}{(1-s)^2+t^2}=\int_{0}^{\infty}e^{(s-1)u}\cos tu\, du.
\end{equation}
This along with similar analysis as above gives (\ref{eq20}). Now it turns out that if we use (\ref{change}) instead of (\ref{useint}) with $j=-1$ when Re $s>1$, then we \textit{do} get the second term on the right-hand side of equation (19) in \cite{riemann}, i.e., $\displaystyle-\frac{1}{4}(4\pi)^{\frac{(s-3)}{2}}\Gamma(s)\zeta(s)\cosh n(1-s)$, as given by Ramanujan. This explains how Ramanujan was erroneously led to his equation.

\textbf{Acknowledgements.} The author wishes to express his gratitude to Professor Bruce C.~Berndt for his constant support, careful reading of this manuscript and for helpful comments. The author would also like to thank Professor Harold G.~Diamond, Boonrod Yuttanan and Jonah Sinick for their help.


\end{document}